\documentclass[11pt]{article}
\usepackage[T1]{fontenc}
\usepackage[utf8]{inputenc}
\usepackage{amsmath,amssymb,amsthm,booktabs,needspace}
\usepackage[margin=1in]{geometry}
\usepackage[colorlinks=true,linkcolor=blue,citecolor=blue,urlcolor=blue]{hyperref}
\newtheorem{theorem}{Theorem}
\newtheorem{lemma}{Lemma}
\newtheorem{corollary}{Corollary}
\newtheorem{proposition}{Proposition}
\theoremstyle{remark}
\newtheorem{remark}{Remark}
\newcommand{\N}{\mathbb Z_{\geq0}}
\newcommand{\SG}{\mathcal G}
\newcommand{\mex}{\operatorname{mex}}
\newcommand{\one}{\mathbf1}

\title{The Sprague--Grundy Function of Slow Exact $k$-Nim\\
with $k+1$ Piles}
\author{%
  \textbf{Vladimir Gurvich}\\[2pt]
  {\small\it National Research University Higher School of Economics (HSE), Moscow, Russia and}%
  \\[2pt]
  {\small\it RUTCOR, Rutgers University, Piscataway, NJ, United States}\\[2pt]
  \texttt{vgurvich@hse.ru}, \texttt{vladimir.gurvich@gmail.com}
  \vspace{10pt}
  \and
  \textbf{Mariya Naumova} (Corresponding author)\\[2pt]
  {\small\it Rutgers Business School, Rutgers University, Piscataway, NJ, United States}\\[2pt]
  \texttt{mnaumova@business.rutgers.edu}
}
\date{}

\begin{document}
\maketitle

\begin{abstract}
We give an explicit formula for the normal-play Sprague--Grundy (SG)
function of slow exact $k$-Nim with $k+1$ piles. A move removes one
stone from each of exactly $k$ nonempty piles. A truncation preserves
every legal move sequence. The SG value depends only on the total of
the truncated coordinates modulo $2k$, their number of odd entries,
and the parities of the two smallest piles. For every $k\geq2$, the
set of the SG values is $\{0,1,2,3\}$. We prove the formula by the mex
characterization, including the reductions needed after boundary
moves. We also prove a unit-shift identity and use it with a known
optimal-move theorem to explain the four-value bound independently
of the explicit formula.
\end{abstract}

\section{Introduction and main result}

The slow versions of exact $k$-Nim and Moore's $k$-Nim were introduced
by Gurvich and Ho~\cite{GH}. For integers $n\geq k\geq1$, a position
is a vector of $n$ nonnegative pile sizes. Two players move alternately.
A slow exact move decreases exactly $k$ positive coordinates by one;
a slow Moore move decreases at least one and at most $k$.
We use normal play, in which a player with no legal move loses;
under mis\`ere play, that player wins.

A position in slow exact $k$-Nim is terminal if and only if fewer than $k$
coordinates are positive, or equivalently, at least $n-k+1$
coordinates are zero. Thus, when $n=k+1$, a position is terminal
if and only if at least two coordinates are zero.
The cases $k=1$ and $k=n$ have fixed play lengths $\sum_i x_i$
and $\min_i x_i$, respectively. We focus on $n=k+1$ with $k\geq2$.

For $k=2$, the normal-play Sprague--Grundy (SG) function and its
mis\`ere analogue were characterized for $n\leq4$ in~\cite{GHHC};
the latter uses the mex recursion with terminal value $1$.
The normal-play $P$-positions for $n=5,6$ were determined
in~\cite{CGKPV}. For $k=2$ and $n=4,5$, remoteness was studied
computationally in~\cite{GMN}, which proposed formulas for the
observed exceptions to the GM-rule but did not prove their completeness.
For $n=k+1$, the normal-play remoteness function and an optimal
move rule were established in~\cite{GMMV}; see also~\cite{GN24}.
Partial results for the mis\`ere version appear in~\cite{GMMN}.

Dufour and Heubach~\cite{DH} characterized the $P$-positions for
$n=k+1$ using reduced positions, the total modulo $2k$, and the
number of odd piles. We use their reduction framework and extend
this characterization to the complete SG function. In particular,
we distinguish the values $1$, $2$, and $3$ and give a direct mex
proof that includes the necessary boundary reductions.

Positions are sorted in nondecreasing order after each move.
For $n=k+1$, a move \emph{omits} the pile it leaves unchanged.
Write $\one=(1,\ldots,1)$. The SG function is
defined by
\[
 \SG_k(x)=\mex\{\SG_k(z):x\longrightarrow z\},
 \qquad \mex\varnothing=0,
\]
where $\mex$ is the least excluded nonnegative integer. Its zero
set is the set of $P$-positions.

\Needspace{30\baselineskip}
\begin{theorem}\label{thm:main}
Let $k\geq2$, let $n=k+1$, and let
\[
 x\in\N^n,\qquad 0\leq x_1\leq\cdots\leq x_n.
\]
For $2\leq j\leq n$, put $P_j=x_1+\cdots+x_j$, and define
\begin{equation}\label{eq:cap}
 L=\min_{2\leq j\leq n}
       \left\lfloor\frac{P_j}{j-1}\right\rfloor,
 \qquad y_i=\min\{x_i,L\}.
\end{equation}
Let
\begin{equation}\label{eq:stats}
 S=\sum_{i=1}^{n}y_i,\qquad
 s=S\bmod 2k,\qquad
 o=\#\{i:y_i\text{ is odd}\}.
\end{equation}
Here $s$ is the least nonnegative residue modulo $2k$.
The two special residues give
\begin{equation}\label{eq:special}
 \boxed{\quad s=k-1\ \Longrightarrow\ \SG_k(x)=0,
 \qquad s=2k-1\ \Longrightarrow\ \SG_k(x)=1.\quad}
\end{equation}
For every other residue, define
\begin{equation}\label{eq:r}
 r=\min\{s,\,2k-2-s\}.
\end{equation}
Then
\begin{equation}\label{eq:formula}
\boxed{
\SG_k(x)=
\begin{cases}
0,&o\leq r,\\
1,&o\geq r+4,\\
1,&o=r+2,\ y_1\text{ even},\\
2,&o=r+2,\ y_1,y_2\text{ odd},\\
3,&o=r+2,\ y_1\text{ odd},\ y_2\text{ even}.
\end{cases}}
\end{equation}
These cases are exhaustive. Moreover, $y_1=x_1$ and $y_2=x_2$,
so the parities of the two smallest piles can be read from the original position.
\end{theorem}

Indeed, $o\equiv s\equiv r\pmod2$ for the nonspecial residues,
so $o-r$ is even. Also $0\leq r\leq k-2$ there.
For completeness, when $k=1$ every move removes one stone and
\[
 \SG_1(x_1,x_2)=(x_1+x_2)\bmod2.
\]
In the proof below we therefore assume $k\geq2$.

\section{Truncation preserves the game tree}
\label{sec:reduction}

We use the reduction framework of Dufour and Heubach~\cite{DH}.
Their Lemma~1 proves preservation of move sequences, Theorem~4
gives the reduced-position criterion, and Section~3.4 describes
initial reduction. For the present family, \eqref{eq:cap} gives
a direct formula for the truncation level. We include the proof
to make the SG argument self-contained.

We call a sorted position $y$ \emph{reduced} when
\begin{equation}\label{eq:reduced}
 \sum_{i=1}^{n}y_i\geq k y_n.
\end{equation}

\Needspace{10\baselineskip}
\begin{lemma}\label{lem:reduction}
The truncation \eqref{eq:cap} preserves every legal move sequence,
with the piles regarded as labeled. In particular,
$\SG_k(x)=\SG_k(y)$. The resulting position is reduced, and its
first two coordinates are unchanged.
\end{lemma}

\begin{proof}
Every move decreases at least $j-1$ of the first $j$ original piles.
Thus a legal sequence of $\ell$ moves satisfies
\[
 (j-1)\ell\leq P_j\qquad(2\leq j\leq n),
\]
and therefore $\ell\leq L$.
No individual pile can lose more than $L$ stones along this sequence.
Replacing a pile larger than $L$ by $L$ cannot invalidate the sequence.
Conversely, a sequence legal from the smaller vector $y$ is legal
from $x$. Thus the rooted game trees are isomorphic under the
correspondence that preserves each sequence of omitted-pile labels.

Let $q$ be the number of original coordinates strictly smaller
than $L$. If $q\geq2$, then
\[
 \sum_i y_i=P_q+(n-q)L
 \geq(q-1)L+(n-q)L=kL\geq k y_n.
\]
For $q=0,1$, the same inequality is immediate.
This proves \eqref{eq:reduced}.

Finally, $P_j\geq(j-1)x_2$ for every $j\geq2$, since
$x_i\geq x_2$ for $i\geq2$ and $x_1\geq0$.
Consequently $L\geq x_2$, proving $y_1=x_1$, $y_2=x_2$.
\end{proof}

\begin{proposition}[Maximum length of play]\label{prop:height}
The integer $L$ in \eqref{eq:cap} is the maximum possible number
of successive legal moves from $x$.
\end{proposition}

\begin{proof}
Lemma~\ref{lem:reduction} shows that every legal play has length
at most $L$. Its proof also establishes the stronger inequality
$\sum_i y_i\geq kL$, while $0\leq y_i\leq L$ for every $i$.
We may therefore choose integers $a_i$ with
$0\leq a_i\leq y_i$ and $\sum_i a_i=kL$.
Put $e_i=L-a_i$. Then $e_i\geq0$ and
\[
 \sum_{i=1}^{k+1}e_i=(k+1)L-kL=L.
\]
Schedule $L$ moves, omitting pile $i$ exactly $e_i$ times.
Pile $i$ is decreased exactly $a_i$ times, so every scheduled
move is legal: before any scheduled decrease, fewer than $a_i$
stones have been removed from that pile. This gives a legal
sequence of length $L$ from $y$, and hence from $x$.
\end{proof}

Every reduced position is unchanged by the truncation. In fact,
if $d=y_n$ and $S\geq kd$, then for every $j\geq2$,
\[
 \sum_{i=1}^{j}y_i
 =S-\sum_{i=j+1}^{n}y_i
 \geq kd-(n-j)d=(j-1)d.
\]
Also, the only terminal reduced position is the zero vector:
if two coordinates vanish and $d>0$, then
$S\leq(k-1)d<kd$.

After each move we may truncate again. Such a move followed by
truncation strictly decreases the total number of stones, so induction
on that total is available.

\Needspace{12\baselineskip}
\section{Moves requiring another truncation}
\label{sec:followers}

The following classification is the specialization of the follower
reductions in \cite[Lemmas~4--6]{DH} needed below. We express it
using the nonnegative slack $t$ and include the details that control
the two smallest coordinates.

Let $y$ be reduced, and put
\begin{equation}\label{eq:slack}
 d=y_n,\qquad S=\sum_i y_i,\qquad t=S-kd\geq0.
\end{equation}

\begin{lemma}\label{lem:followers}
A move omitting a nonmaximal pile is already reduced afterward.
For a legal move omitting a maximal pile, the reduced follower is:
\begin{enumerate}
\item the ordinary follower if $t\geq k$;
\item $y-\one$ if $1\leq t\leq k-1$;
\item if $t=0$, the vector obtained by decreasing all nonmaximal
      piles by one and all maximal piles by two.
\end{enumerate}
In the last case, if $\alpha$ is the number of maximal piles,
then $1\leq\alpha\leq k-1$.
In both cases requiring truncation, the two smallest pile
parities flip.
\end{lemma}

\begin{proof}
If a nonmaximum is omitted, all maximal piles decrease to $d-1$,
and the ordinary follower has total $S-k$. Its reduced slack is
$(S-k)-k(d-1)=t\geq0$.
If a maximum is omitted, its size remains $d$ and the ordinary
slack is $t-k$. This proves the first assertion and case~1.

We use the following observation to identify the exact truncation
level in the remaining two cases. Let $w$ be the sorted ordinary
follower, and let $L(w)$ denote its truncation level as defined
in \eqref{eq:cap}. Suppose that capping $w$ at a nonnegative
integer $c$ produces a reduced vector
\[
 v_i=\min\{w_i,c\},\qquad v_n=c.
\]
Then, for every $2\leq j\leq n$,
\[
 \sum_{i=1}^{j}w_i
 \geq \sum_{i=1}^{j}v_i
 = \sum_{i=1}^{n}v_i-\sum_{i=j+1}^{n}v_i
 \geq kc-(n-j)c
 = (j-1)c,
\]
where we used the reducedness of $v$ and $n=k+1$.
Consequently,
\[
 L(w)=\min_{2\leq j\leq n}
 \left\lfloor\frac{\sum_{i=1}^{j}w_i}{j-1}\right\rfloor
 \geq c.
\]
Thus, if a prefix bound also gives $L(w)\leq c$, then
$L(w)=c$, proving that $c$ is the exact truncation level.
We apply this observation below with $c=d-1$ and $c=d-2$,
respectively.

For $1\leq t\leq k-1$, capping the omitted maximum at $d-1$
produces $y-\one$, whose slack is
\[
 S-(k+1)-k(d-1)=t-1\geq0.
\]
The first $k$ piles of the ordinary follower have sum
$(k-1)d+t-k$, so their prefix bound forces a cap of at most
\[
 \left\lfloor\frac{(k-1)d+t-k}{k-1}\right\rfloor=d-1.
\]
The capped vector is reduced and has maximum $d-1$, so the preceding
observation proves that $d-1$ is the exact cap.

Suppose now that $t=0$. If $\alpha\geq k$, then $S=kd$
would force a remaining pile to be zero, or the entire position
to be zero. Omitting a maximum would then be illegal.
Thus $\alpha\leq k-1$. Moreover, legality implies that every
pile is positive, so $d\geq2$ when $t=0$.
The same prefix bound forces a cap of at most $d-2$.
Capping there decreases all $\alpha$ original maxima by two
and every other pile by one. Its slack is
\[
 S-(k+1)-\alpha-k(d-2)=k-1-\alpha\geq0,
\]
so the capped vector is reduced and has maximum $d-2$.
The same observation proves that $d-2$ is the exact cap.

In case~2 every pile decreases by one.
In case~3 at least two piles are nonmaximal, so the first two
decrease by one as well. Their order is preserved after capping.
\end{proof}

\section{The \texorpdfstring{$\delta$}{delta}-transformation}

We call a residue \emph{regular} if
$s\notin\{k-1,2k-1\}$. For a reduced position with regular
residue, define the integer
\begin{equation}\label{eq:delta}
 \delta=\frac{o-r}{2},\qquad r=\min\{s,2k-2-s\}.
\end{equation}

\begin{lemma}\label{lem:delta}
Let a move omit a pile of parity $p\in\{0,1\}$, with $p=1$
meaning odd. If the current residue is regular, the reduced
follower also has regular residue and satisfies
\begin{equation}\label{eq:delta-move}
 \boxed{\delta'=1-\delta+p.}
\end{equation}
\end{lemma}

\begin{proof}
If no truncation is required, then
\[
 s'=(s-k)\bmod2k,\qquad o'=k-o+2p.
\]
Regular residues remain regular, and
$r'=k-2-r$. Hence
\[
 \delta'=\frac{k-o+2p-(k-2-r)}2=1-\delta+p.
\]

It remains to check the moves in cases~2 and~3 of
Lemma~\ref{lem:followers}. Here the omitted maximum has parity
$p=d\bmod2$. Since $0\leq t<k$,
\[
 s=\begin{cases}t,&p=0,\\k+t,&p=1.\end{cases}
\]
For a regular current residue, case~2 has $1\leq t\leq k-2$.
The reduced follower's statistics are as follows:
\begin{equation}\label{eq:boundary-table}
\begin{array}{c|c|c|c|c}
 \text{Condition}&p&s'&o'&r'\\ \hline
 1\leq t\leq k-2&0&k+t-1&k+1-o&k-1-t\\
 1\leq t\leq k-2&1&t-1&k+1-o&t-1\\
 t=0&0&k-1-\alpha&k+1-o-\alpha&k-1-\alpha\\
 t=0&1&2k-1-\alpha&k+1-o+\alpha&\alpha-1
\end{array}
\end{equation}
The first two rows follow from $y'=y-\one$.
For the last two, maximal piles decrease by two, so their parities
stay unchanged, while the other $k+1-\alpha$ parities flip.
Since $1\leq\alpha\leq k-1$, all displayed residues are regular.
Substituting the entries into $(o'-r')/2$ gives
$1-\delta$ when $p=0$, and $2-\delta$ when $p=1$.
\end{proof}

\section{Proof of Theorem~\ref{thm:main}}

Let $F(y)$ denote the proposed value for a reduced position.
We prove $F=\SG_k$ by strong induction on $S=\sum_i y_i$.
The terminal reduced vector is zero. For $k\geq2$, it has
$s=o=r=0$ and therefore $F=0$.
All reduced followers have smaller total, by
Lemma~\ref{lem:reduction}.

\subsection{Regular residues with
\texorpdfstring{$\delta\leq0$ or $\delta\geq2$}{delta <= 0 or delta >= 2}}

If $\delta\leq0$, equation~\eqref{eq:delta-move} gives
$\delta'\geq1$ for every follower. Every follower therefore
has positive proposed value, so $\mex=0=F(y)$.

If $\delta\geq2$, the same equation gives $\delta'\leq0$
for every follower. All followers have value zero.
The position is nonterminal, since $o\geq r+4\geq4$.
Consequently $\mex\{0\}=1=F(y)$.

\subsection{Regular residues with \texorpdfstring{$\delta=1$}{delta = 1}}

Here $o=r+2$ and $0\leq r\leq k-2$, so
\begin{equation}\label{eq:odd-range}
 2\leq o\leq k.
\end{equation}
There are at least two odd piles and at least one even pile.
Omitting an even pile gives $\delta'=0$, hence value zero.
Omitting an odd pile gives $\delta'=1$.

\paragraph{The smallest pile is even.}
Omitting it is legal: a nonzero reduced position has at most
one zero pile. This move gives value zero.
Every legal omission of an odd pile decreases a smallest pile from
even to odd. The follower stays in the case $\delta'=1$ and
therefore has proposed value $2$ or $3$.
This remains true for a move requiring truncation by
Lemma~\ref{lem:followers}. No follower has value $1$, so
the mex is $1$.

\paragraph{The two smallest piles are odd.}
All piles are positive. Omitting an even pile gives value zero.
Omitting an odd pile other than the first leaves a smallest
pile to be reduced from odd to even, giving value $1$.
Such a pile exists by \eqref{eq:odd-range}. If this move requires
truncation, Lemma~\ref{lem:followers} ensures that the smallest
entry still decreases by one.

Every odd omission has value $1$ unless it retains a unique
smallest pile. In that latter case, the smallest pile stays odd,
and the second-smallest decreases from odd to even; the follower
has value $3$. If the minimum is repeated, omitting one copy
still reduces another copy to an even minimum, giving value $1$.
Thus values $0$ and $1$ occur, and value $2$ does not.
The mex is $2$.

\paragraph{The smallest pile is odd and the second is even.}
The minimum is unique. Omitting it leaves the minimum odd and
decreases the second-smallest to odd, so the follower has value $2$.
Omitting another odd pile gives an even minimum and value $1$;
this remains true after any required truncation by
Lemma~\ref{lem:followers}.
Omitting an even pile gives value zero.
These are all possible follower values. Hence the mex is $3$.

\subsection{The special residue \texorpdfstring{$s=k-1$}{s = k - 1}}

Every move not requiring truncation has follower residue $2k-1$
and therefore value $1$.
If a move requires truncation, then necessarily
\[
 t=k-1,\qquad d\text{ even},
\]
and its reduced follower is $y-\one$. This follower has
\[
 s'=2k-2,\qquad r'=0,\qquad o'=k+1-o.
\]
The original position has at least one even pile, namely a maximum.
Also $o\equiv k-1\equiv k+1\pmod2$, so the number of even piles
is positive and even. It is therefore at least two.
Thus $o'\geq2$, and the follower has positive proposed value.
No move reaches value zero, proving $F(y)=0$.

\subsection{The special residue \texorpdfstring{$s=2k-1$}{s = 2k - 1}}

Every move not requiring truncation has follower residue $k-1$
and value zero.
If a move requires truncation, then necessarily
\[
 t=k-1,\qquad d\text{ odd},
\]
and its reduced follower is $y-\one$, with
\[
 s'=k-2,\qquad r'=k-2,\qquad o'=k+1-o.
\]
Since $s$ is odd, $o$ is odd. If $o\geq3$, then
$o'\leq k-2=r'$, giving value zero.
If $o=1$, the unique odd pile is the unique maximum.
The two smallest original piles are even, so the two smallest
follower piles are odd. Also $o'=k=r'+2$, giving follower value $2$.

In this last situation all piles are positive, because the move
omitting a maximum was legal. An even nonmaximal pile can instead
be omitted. That move needs no truncation and has follower value zero.
Thus a zero-valued follower always exists, while no follower has
value $1$. The mex is $1$.

\qed

\section{Consequences and examples}
\label{sec:consequences}

\begin{corollary}\label{cor:range}
For every $k\geq2$, the range of $\SG_k$ is exactly
$\{0,1,2,3\}$. If the smallest pile is even, its SG value
belongs to $\{0,1\}$.
\end{corollary}

\begin{proof}
The upper bound and the assertion about an even minimum follow
immediately from the formula. The following families realize
all four values for every $k\geq2$:
\[
\begin{array}{c|c}
\text{Position of length $k+1$}&\SG_k\\ \hline
(0,\ldots,0)&0\\
(1,\ldots,1)&1\\
(1,1,\underbrace{2,\ldots,2}_{k-1})&2\\
(1,2,\underbrace{3,\ldots,3}_{k-1})&3
\end{array}
\]
The last two positions are reduced, with $(s,o)=(0,2)$
and $(s,o)=(k,k)$, respectively.
\end{proof}

The zero set recovers \cite[Theorem~6]{DH}. After truncation,
all $P$-positions are given by
\[
 s=k-1\quad\text{or}\quad
 \bigl(s\notin\{k-1,2k-1\}\text{ and }o\leq r\bigr).
\]
The only positions with values $2$ or $3$ satisfy
$o=r+2$ and have odd minimum. The parity of the second-smallest pile
distinguishes those two values.

\paragraph{A zero minimum.}
If $x_1=0$, every move must omit the first pile. Play therefore
lasts exactly $x_2$ moves, including zero moves when $x_2=0$, and
\begin{equation}\label{eq:zero-minimum}
 \SG_k(0,x_2,\ldots,x_{k+1})=x_2\bmod2.
\end{equation}
The prefix with $j=2$ gives $L\leq x_2$, while
Lemma~\ref{lem:reduction} gives $L\geq x_2$. Thus
$L=x_2$ and $y=(0,x_2,\ldots,x_2)$, so the main formula
reproduces this forced-play case.

\paragraph{Four piles and further examples.}
For $k=3$, Theorem~\ref{thm:main} gives the following formula
for slow exact $(4,3)$-Nim:
\[
\begin{array}{c|c}
s&\text{Rule}\\ \hline
0,4&r=0\\
1,3&r=1\\
2&\SG_3=0\\
5&\SG_3=1.
\end{array}
\]
Use \eqref{eq:formula} in the first two rows.
Table~\ref{tab:examples} gives examples for several values of $k$.

\begin{table}[htbp]
\centering
\caption{Examples of truncation and SG evaluation.}
\label{tab:examples}
\begin{tabular}{@{}cllcc@{}}
\toprule
$k$&$x$&Reduced $y$&$(s,o)$&$\SG_k(x)$\\
\midrule
3&$(1,2,5,6)$&$(1,2,3,3)$&$(3,3)$&3\\
4&$(1,1,2,2,2)$&$(1,1,2,2,2)$&$(0,2)$&2\\
4&$(1,2,3,3,100)$&$(1,2,3,3,3)$&$(4,4)$&3\\
5&$(1,1,2,2,2,2)$&$(1,1,2,2,2,2)$&$(0,2)$&2\\
5&$(1,2,3,3,3,3)$&$(1,2,3,3,3,3)$&$(5,5)$&3\\
\bottomrule
\end{tabular}
\end{table}

\begin{remark}[A correction to an earlier computation]\label{rem:old-computation}
Remark~16 of the 2015 preprint~\cite{GH} reports SG value $5$
for the positions $(1,2,3,3,m)$ of slow exact $(5,4)$-Nim when
$3\leq m\leq7$. In fact, the value is $3$ for every $m\geq3$.
This can be checked directly, without Theorem~\ref{thm:main}.
Omitting the pile of size $1$ gives $(1,1,2,2,m-1)$.
Its followers have SG values $0$ and $1$: decreasing both piles
of size $1$ gives a terminal position, while omitting one of them
gives a forced path of length one. Thus its SG value is $2$.
Omitting the pile of size $2$ gives $(0,2,2,2,m-1)$, of value
$0$ by \eqref{eq:zero-minimum}. Every other omission gives
a position with a zero minimum and second-smallest entry $1$,
of value $1$. Hence the option values are exactly $\{0,1,2\}$,
and the original position has SG value $3$.
\end{remark}

There is also a uniform version of the residue-and-gap description.
For reduced $y$, put $g_i=y_{i+1}-y_i$ and $t=S-ky_n$.
Then
\begin{equation}\label{eq:gap-coordinates}
 y_1=t+\sum_{i=1}^{k}(i-1)g_i,\qquad
 S=(k+1)t+k\sum_{i=1}^{k}i g_i.
\end{equation}
Consequently the SG value depends only on $t\bmod2k$ and
the parities of all the gaps. The theorem evaluates this dependence
without constructing a table with exponentially many parity patterns.

\Needspace{14\baselineskip}
\section{A structural explanation of the four-value bound}
\label{sec:four-values-explanation}

The preceding mex proof is self-contained. We now give a separate
explanation of the bound in Corollary~\ref{cor:range}, using
a unit-shift identity and the known normal-play optimal-move
theorem. Throughout this section, $n=k+1$ and $k\geq2$.

\subsection{A unit-shift identity}

For a sorted position, define
\[
 T(x_1,\ldots,x_n)=(x_1,x_2+1,\ldots,x_n+1).
\]
Thus $T$ adds one stone to every pile except a smallest one.
Chikin et al.~\cite[Theorem~3]{CGKPV} proved that adding
$(0,2,\ldots,2)$ preserves SG values in this family. The following
identity refines that normal-play invariance.

\begin{proposition}[Unit shift]\label{prop:unit-shift}
For every sorted position $x\in\N^{k+1}$,
\begin{equation}\label{eq:explanation-unit-shift}
 \boxed{\SG_k(Tx)=\SG_k(x)\oplus1,}
\end{equation}
where $\oplus$ denotes bitwise exclusive-or.
\end{proposition}

\begin{proof}
Write $\operatorname{Opt}(x)$ for the set of sorted legal
followers of $x$, with duplicates ignored. The move rule gives
\begin{equation}\label{eq:explanation-option-shift}
 \operatorname{Opt}(Tx)
 =\{x\}\cup\{Tz:z\in\operatorname{Opt}(x)\}.
\end{equation}
Indeed, leaving the first pile of $Tx$ unchanged gives $x$.
For $i>1$, leaving pile $i$ unchanged gives $Tz$, where $z$ is
the follower of $x$ obtained by leaving that same pile unchanged;
such moves are legal exactly when $x_1>0$.
Conversely, if $z$ is obtained from $x$ by leaving the first pile
unchanged and $x_1<x_2$, then $Tz=x$.
If $x_1=x_2>0$, the same sorted follower $z$ is obtained by
leaving a different smallest pile unchanged, so its image is
already included. If $x_1=x_2=0$, that move is illegal.
This also explains why sorting and ignoring duplicates are
needed in \eqref{eq:explanation-option-shift}.

For any finite set $U\subseteq\N$, put $m=\mex U$. Then
\begin{equation}\label{eq:explanation-mex-shift}
 \mex\bigl(\{m\}\cup\{u\oplus1:u\in U\}\bigr)=m\oplus1.
\end{equation}
If $m$ is even, the set on the left contains $0,\ldots,m$ and
misses $m+1$; if $m$ is odd, it contains every nonnegative
integer below $m-1$ and misses $m-1$.
Induction on $\sum_i x_i$, using
\eqref{eq:explanation-option-shift} and
\eqref{eq:explanation-mex-shift}, proves
\eqref{eq:explanation-unit-shift}.
\end{proof}

In particular, $T$ exchanges the pairs of SG values $0,1$;
$2,3$; $4,5$; and so on. Applying the proposition twice gives
$\SG_k(T^2x)=\SG_k(x)$, recovering the normal-play case of
the previously known double-shift invariance.

\Needspace{16\baselineskip}
\subsection{Even minima and the number of option values}

\begin{corollary}[Even minimum]\label{cor:even-minimum}
For every sorted position $x$,
\begin{equation}\label{eq:explanation-even-minimum}
 \boxed{x_1\text{ even}\quad\Longrightarrow\quad
        \SG_k(x)\in\{0,1\}.}
\end{equation}
\end{corollary}

\begin{proof}
Now suppose $x_1$ is even. In $Tx$, this pile is uniquely
smallest and even, including when $x_1=0$.
The normal-play optimal-move theorem \cite[Theorem~2]{GMMV} prescribes leaving a
smallest even pile unchanged and decreasing every other pile.
Thus its prescribed move from $Tx$ is exactly $Tx\longrightarrow x$.
If $Tx$ is an $N$-position, this optimal move reaches a
$P$-position. If $Tx$ is a $P$-position, every legal move reaches
an $N$-position. Consequently exactly one of $Tx$ and $x$ is a
$P$-position, so exactly one of their SG values is zero.
By \eqref{eq:explanation-unit-shift}, their SG values differ by
xor with $1$, which forces
\[
 \{\SG_k(Tx),\SG_k(x)\}=\{0,1\}.
\]
\end{proof}

This proof uses Proposition~\ref{prop:unit-shift} and
\cite[Theorem~2]{GMMV}, but does not use the explicit SG formula.

\paragraph{At most one exceptional move.}
Suppose that the smallest pile has odd size $a$.
Every move that decreases a smallest pile produces a position
whose minimum is $a-1$, which is even. By
\eqref{eq:explanation-even-minimum}, such a follower has SG value
$0$ or $1$.
The only way to avoid this is to leave a unique smallest pile
untouched. There is at most one such move. If its follower has
SG value $h$, then
\[
 \{\SG_k(z):x\longrightarrow z\}\subseteq\{0,1,h\}.
\]
If there is no exceptional move, the option values are a subset
of $\{0,1\}$. In either case, there are at most three distinct
option values, and therefore $\SG_k(x)\leq3$.
Together with Corollary~\ref{cor:even-minimum}, this proves the
bound for every position.

Equivalently, an SG value of $4$ would require followers of all
four values $0,1,2,3$. Both values $2$ and $3$ would have to come
from the single exceptional move, which is impossible. Thus
increasing the number of piles can increase the number of moves
without increasing the SG bound.
If the odd minimum occurs at least twice, every move decreases
at least one of its copies. There is then no exceptional move,
giving the stronger conclusion $\SG_k(x)\leq2$.

\paragraph{Why the value three occurs.}
For $(n,k)=(4,3)$, the distinct followers of $(1,2,3,3)$
are listed in Table~\ref{tab:four-pile-options}.
\begin{table}[htbp]
\centering
\caption{The option values that give $\SG_3(1,2,3,3)=3$.}
\label{tab:four-pile-options}
\begin{tabular}{@{}lcc@{}}
\toprule
Pile left unchanged&Sorted follower&$\SG_3$\\
\midrule
The pile of size $1$&$(1,1,2,2)$&2\\
The pile of size $2$&$(0,2,2,2)$&0\\
Either pile of size $3$&$(0,1,2,3)$&1\\
\bottomrule
\end{tabular}
\end{table}
Hence $\SG_3(1,2,3,3)=\mex\{0,1,2\}=3$.
These follower values can be checked directly. With exactly one
zero pile, the move is forced: that pile must be left unchanged
until a second pile becomes zero. Thus $(0,2,2,2)$ has a forced
path of length two and SG value $0$, while $(0,1,2,3)$ has a
forced path of length one and SG value $1$.
The two distinct followers of $(1,1,2,2)$ are $(0,1,1,1)$ and
$(0,0,1,2)$, of values $1$ and $0$, respectively, so its value
is $2$. The families in the proof of Corollary~\ref{cor:range}
show that all four values occur for every $k\geq2$.

\Needspace{16\baselineskip}
\section{Complexity and implementation}
\label{sec:implementation}

The input is an explicit list of $n=k+1$ nonnegative integers
written in binary. For a sorted input, the prefix sums, truncation, total, odd count,
and final case distinction require $O(n)$ arithmetic operations.
Sorting adds $O(n\log n)$ comparisons.
If each input coordinate has at most $B$ bits, all intermediate
integers have $O(B+\log n)$ bits. Hence the full SG function is
polynomial-time computable even when $n$, and therefore $k$, varies.
For fixed $n$, all divisors and moduli are constants, and the
bit complexity is $O(B)$ in the usual bit model.

For implementation, we use the following compact form.
Define
$r=\min\{s,2k-2-s\}$ for all residues, including $r=-1$
when $s=2k-1$:
\begin{equation}\label{eq:compact}
\SG_k(x)=
\begin{cases}
0,&s=k-1\text{ or }o\leq r,\\
3-(x_2\bmod2),&s\notin\{k-1,2k-1\},\ o=r+2,
                 \ x_1\text{ odd},\\
1,&\text{otherwise}.
\end{cases}
\end{equation}

\section{Exact computational verification}
\label{sec:verification}

A direct mex computation over the original legal moves agreed with
the formula on every sorted position in Table~\ref{tab:verification}.
For $n$ piles with sizes at most $M$, the number of sorted positions
is $\binom{M+n}{n}$. The first eleven rows contain $2\,174\,629$
positions; the last two rows extend the check to a total of
$2\,213\,389$ positions. The recursive SG relation was evaluated
in lexicographic order with stored follower values, without using
truncation in the original-game computation.

\begin{table}[htbp]
\centering
\caption{Exhaustive verification using the original legal moves;
$k=n-1$ in every row.}
\label{tab:verification}
\begin{tabular}{@{}rrr@{}}
\toprule
$n$&Maximum pile size $M$&Number of positions\\
\midrule
2&40&861\\
3&65&50\,116\\
4&40&135\,751\\
5&35&658\,008\\
6&24&593\,775\\
7&17&346\,104\\
8&13&203\,490\\
9&10&92\,378\\
10&8&43\,758\\
11&7&31\,824\\
12&6&18\,564\\
\midrule
13&6&27\,132\\
14&5&11\,628\\
\midrule
\multicolumn{2}{r}{Total}&2\,213\,389\\
\bottomrule
\end{tabular}
\end{table}

\section*{Statement about AI usage} 

ChatGPT was used during the preparation of the manuscript to assist
with coding, exploratory computations, the development of mathematical
arguments (in particular, Theorem 1), and partial editing. The authors take full responsibility for all
content in this paper.

\section*{Acknowledgments} 
This paper is an output of the research project (HSE-BR-2025-024) 
implemented as part of the Basic Research Program at HSE University.

\section*{Funding}
The authors declare that no financial support was received for the
research, authorship, or publication of this article.

\section*{Disclosure of interest}
The authors report no competing interests to declare.

\end{document}